\documentclass[11pt]{article}

\usepackage[colorlinks=true,pagebackref,hyperindex]{hyperref}
\usepackage{mathrsfs} 
\usepackage[all]{xy}
\usepackage{color}
\usepackage{amsfonts}
\usepackage{amscd}

\definecolor{darkgreen}{rgb}{0.0, 0.7, 0.0}

\definecolor{cyan}{cmyk}{1,0,0,0}

\newtheorem{tm}{Theorem}[subsection]
\newtheorem{lm}[tm]{Lemma}
\newtheorem{pr}[tm]{Proposition}
\newtheorem{rmk}[tm]{Remark}
\newtheorem{cor}[tm]{Corollary}

\newtheorem{??}[tm]{Question}

\newtheorem{num}[tm]{}

\newcommand\blacksquare{{\hspace*{\fill} $\fbox{}$}} 
\newcommand{\ben}{\begin{enumerate}}
\newcommand{\een}{\end{enumerate}}
\newcommand{\bit}{\begin{itemize}}
\newcommand{\eit}{\end{itemize}}
\newcommand{\beq}{\begin{equation}}
\newcommand{\eeq}{\end{equation}}
\newcommand{\la}{\label}
\newcommand{\n}{\noindent}
\newcommand\ci{\cite}

\font\tenmsb=msbm10
\font\sevenmsb=msbm7
\font\fivemsb=msbm5
\newfam\msbfam
\textfont\msbfam=\tenmsb
\scriptfont\msbfam=\sevenmsb
\scriptscriptfont\msbfam=\fivemsb
\def\Bbb#1{{\fam\msbfam #1}}
\font\teneufm=eufm10
\font\seveneufm=eufm7
\font\fiveeufm=eufm5
\newfam\eufmfam
\textfont\eufmfam=\teneufm
\scriptfont\eufmfam=\seveneufm
\scriptscriptfont\eufmfam=\fiveeufm
\def\frak#1{{\fam\eufmfam\relax#1}}

\newcommand{\ke}{ \hbox{\rm Ker} }
\newcommand{\lorw}{\longrightarrow}

\newcommand\rat{{\Bbb Q}}

\newcommand\oql{\overline{\Bbb Q}_\ell}
\newcommand\comp{{\Bbb C}}
\newcommand\real{{\Bbb R}}
\newcommand\zed{{\Bbb Z}}

\newcommand\s{\sigma}

\newcommand{\ov}[1]{\overline{#1}}

\newcommand{\m}[1]{\mathcal{#1}}

\newcommand{\bb}[1]{\Bbb{#1}}

\newcommand{\z}[1]{{_0{#1}}}

\newcommand{\om}{\Omega}

\begin{document}

\title{Proper toric maps over finite fields}

\author{
Mark Andrea A.  de Cataldo\thanks{
Partially supported by N.S.F. grant DMS-1301761 and by a grant from the Simons Foundation (\#296737 to Mark de Cataldo)}
}

\date{2014}

\maketitle

\begin{abstract}
We determine 
 a strong form of the  decomposition theorem
 for proper toric maps over finite fields.
 \end{abstract}

\tableofcontents

\section{Introduction, notation, basic toric geometry, statements}\la{int}

\subsection{Introduction}\la{intro}
For a general proper map
of varieties over a finite field, the  decomposition theorem in \ci{BBD}
predicts that  the direct image of the intersection cohomology complex 
becomes semisimple after passage to an algebraic closure. In this paper, we prove Theorem \ref{tmb} which establishes that, for proper toric maps of toric varieties, the above semisimplicity already occurs over the finite field. The simple direct summands
are described explicitly. Recall that the semisimplicity is not known
 even for the cohomology of smooth projective varieties over finite fields.

If the  proper toric fibration  is surjective
with connected fibers, 
then the direct summands
appearing in the decomposition theorem are Tate-twisted intersection
cohomology complexes of closures of orbits on the target with constant coefficients;
see Theorem \ref{tmb}.c.
For a general proper toric map $f,$ we first Stein factorize the map $f = h \circ g,$
we apply the above
to the proper toric fibration $g,$  we push forward each direct summand 
via the finite toric  map $h$ and examine the result; see Theorem \ref{tmb}.a
and Lemma \ref{rt5m}.

In order to carry out what above, we first prove Theorem \ref{tmb}.b, which here we state in the special case when  the domain is smooth: 
the cohomology of the fibers over closed points is trivial in odd degree, and in even 
degree is pure with eigenvalues of Frobenius given by a suitable power
of the cardinality of the finite field. The triviality in odd degree, for example,
implies that, in the context of proper toric maps,  a plane cubic with a node  may not appear as  a fiber of a proper toric map (of course, it may appear as the image of one).

This paper is a companion paper to \ci{dMM}, where the precise form of the
decomposition theorem for proper toric varieties over $\comp,$ as well as an analogue of the  purity
statement above,   are proved with different methods. The main purpose
of \ci{dMM} is to then  introduce  a topological/combinatorial  invariant 
of proper toric fibrations that detects, for example,  whether a given
orbit contributes a direct summand to the decomposition theorem. This turns out to
be related to seemingly subtle combinatorial positivity questions.

The purpose of this paper is also to offer a sample computation  in the context
of $\oql$-adic cohomology over finite fields in a manner which we hope is accessible
to the non-expert. There are two main  differences with the situation over $\comp:$
1) as mentioned above, in general the  conclusion  of decomposition theorem over an algebraic closure of
a finite field does not seem to hold in its full-strength 
over a finite field (see \S\ref{mn}); 2) even if toric maps
are defined over $\zed,$ the local systems appearing as coefficients
in the decomposition theorem depend on the characteristic of the ground field
(see Remark \ref{rmkchar}).

\subsection{Notation}\la{notz}
This paper deals with   proper toric maps
of toric varieties over a finite field. Our main references are \ci{BBD, Ful93}.

Toric varieties and toric maps
are  defined over
the ring of integers, hence over any ground field.
We view a toric variety
 $Z=Z(\Delta)$  as the one associated with a fan $\Delta$ in a lattice
$N \cong \zed^n$  so that $\dim Z =n.$  If helpful, we add subscripts: $\Delta_Z,$ etc. We view $\Delta$ also as a poset:
$\tau \leq \s$ iff  the cone $\tau$ is a face of the cone  $\s$ iff 
$V(\tau) \supseteq V(\s)$ (reversed inclusion for the closures of the orbits $O(\tau)$ and $O(\s)$).
The support of the fan $\Delta$  in $N_{\real}$ is denoted
by $|\Delta |.$ 
The $n$-dimensional torus $T\subseteq Z$ acts on $Z$ with smooth action map ${\rm act}: Z\times T \to Z.$
 Each orbit $O(\s)$  carries a distinguished point
$z_\s$ which is rational over the prime subfield of the ground field.
 We denote
by $\om$ the resulting partition
 $Z= \coprod_{\s} O(\s)$ into locally closed smooth  subvarieties.

A toric map is a   toric map  $f: X\to Y$ of toric varieties, i.e.
the one  associated with    a map $f_N:
N_X  \to N_Y$ of the lattices  with the following property:
every cone $\s \in \Delta_X$  has image contained inside a cone of $\Delta_Y.$
This gives rise to the map  of posets $f_{\Delta}: \Delta_X \to \Delta_Y,$ sending  $\s \mapsto \ov{\s}:=$ 
the smallest cone in $\Delta_Y$ containing the image $f_N(\s).$ If $\s \in \Delta_X,$ then $f(x_\s) = y_{\ov{\s}}.$ 

A toric map is proper iff $f_{N_\real}^{-1} (| \Delta_Y |)  = | \Delta_X  |.$ 
A proper  toric fibration is a proper toric map $f:X \to Y$ such that $f_* \m{O}_X = \m{O}_Y.$ If $f$ is proper toric, then $f$ is a proper toric fibration 
iff  $f_N$  is surjective.
A proper toric fibration  is surjective with connected fibers. A proper toric map
which is surjective and with connected fibers is not necessarily a proper toric fibration
(e.g. Frobenius).

Unless mentioned otherwise, 
we work with schemes (separated  and  of finite type)  over
a  finite field
 $\z{\bb{F}},$ of which we fix an algebraic closure $\bb{F}.$

 We denote  schemes over $\z{\bb{F}}$  by using
 the pre-fix $\z{-},$  which we remove after pulling-back to $\bb{F}:$  e.g. if $\z{f}: \z{X} \to \z{Y}$ is an $\z{\bb{F}}$-map, then
 we can pull it back to the $\bb{F}$-map $f: X \to Y$; similarly, for the complexes below. A standard notation is $X_0$ etc.; we depart from it for graphical
 reasons.

 We work with the ``derived" category  of mixed complexes
 $D^b_m(\z{X}, \oql),$ endowed with
 the middle perversity $t$-structure \ci{BBD}, p.126, p.101, p.71,   whose elements we simply call  complexes.

 By graded vector space $M^*$, we mean $\zed$-graded: $M^*= 
 \oplus_{j\in \zed} M^j =M^{\rm even}
\oplus M^{\rm odd}.$ We say that $M^*$ is even if $M^{\rm odd} = \{0\}.$
 
 Given a finite extension  $\z{\bb{F}} \subseteq  {_1{\bb{F}}} \subseteq 
 \bb{F},$ we have the open inclusion of Galois  (pro-finite) groups
 $_1G := {\rm Gal}  (\bb{F}/{_1{\bb{F}}}) \subseteq 
 \z{G} := {\rm Gal}  (\bb{F}/{\z{\bb{F}}}).$
 If $\z{C}$ is a complex on $\z{X},$ then the $\zed$-graded object
 $H^*(X,C)$ is a continuous\footnote{continuity will not be mentioned further}
 $\z{G}$-module. If $x \in \z{X} (\bb{F})$ is a closed point
 with residue field a finite extension $_1{\bb{F}},$ then the
 graded object $\m{H}^*(C)_x$ is a  $_1G$-module.
 The weight-like properties of the cohomology groups and stalks we consider
are well-defined independently of the
finite field extension one works with; see [BBD], 5.1.12.
Instead of insisting on $\z{G},$ $_1G,$ etc.,   by abuse of notation,
we simply talk about $G$-modules and their weights.
  
  Given a 
variety  $\z{X},$ we  have the associated  
 shifted intersection complex $\m{I}_{\z{X}};$
 if $\z{X}$ is smooth, then $\m{I}_{\z{X}} = {\oql}_{\z{X}}.$ 
 For convenience, we also use the intersection complex  $I_{\z{X}}:= 
 \m{I}_{\z{X}} [\dim {\z{X}}],$ 
 which is a perverse sheaf on $\z{X}.$ 
 The intersection cohomology
 groups  of a variety $\z{X}$ are   $I\!H^*(X):= H^*(X, \m{I}_X).$

We have the notions of  even, mixed,  and pure $G$-module $M^*:$ e.g.
$M^*$ is said to be pure of weight $w$ if $M^j$ is pure of weight
$w+j$ for every $j.$

A graded $G$-module $M^*$  is said to be Tate if
it is even and each $M^{2k} \cong \oql^{\oplus s_k} (-k)$,  for some
$s_k \in \zed^{\geq 0}$. 

A complex  $\z{C}$ on $\z{X}$ is said to be
punctually pure of weight $w$ if the graded $G$-modules $\m{H}^*(C)_x$ are
pure of weight $w$ for every closed point $x \in \z{X} (\bb{F}).$ 
Similarly, we have the notion of $\z{C}$ being even, and of $\z{C}$
being Tate. In particular,
we have the notion of $\z{C}$ being pure, punctually pure, even and Tate;
e.g. see Theorem \ref{tmb}.

By a result of O. Gabber, 
the intersection complex $\m{I}_{\z{X}}$ of a variety  is pure of weight zero.
 The 
Tate-shifted $\m{I}_{\z{X}}(-k)$ is pure of weight $2k,$ and $I_X$ is pure of weight
$\dim X.$  By a result of P. Deligne,   if $\z{f}: \z{X} \to \z{Y}$ is  a proper map of varieties,
then $R \,\z{f}_*  \m{I}_{\z{X}}$ is pure of weight zero.
However, in  general,   $\m{I}_{\z{X}}$ and  $R\, {\z{f}_*}  \m{I}_{\z{X}}$
are neither punctually pure, even, nor Tate.

A complex on a toric variety is said to be $\om$-constructible
if its restriction to each orbit has lisse\footnote{the notion of lisse sheaf is the
$\oql$-adic analogue of a locally constant sheaf}  cohomology sheaves. A 
skyscraper constant sheaf at the origin of the affine line is $\om$-constructible,
whereas one at the point $1$ is not. The intersection complex 
of a toric variety $\z{Z}$  is $\om$-constructible.
The direct image complex $R\,\z{f}_* \m{I}_{\z{X}}$ via a  proper toric map $\z{f},$
may fail to be $\om$-constructible, e.g. the second  closed embedding above
(the first one is not a toric map, according to the definitions).

Given a toric map $\z{f}: \z{X} \to \z{Y}$ 
 define the cohomology graded sheaf on $Y$
$R^*: = \oplus_j R^j f_* \m{I}_X.$  Denote  the  restriction of $R^*$ to an orbit
$O(\s)\subseteq Y$ 
by $R^*_\s,$ and its stalk at a closed point  $y \in \z{Y} (\bb{F})$  by $R^*_y$ 
--this is a graded $G$-module--.  If $\z{f}$ is proper, then proper base change yields
$R^*_y = H^* (f^{-1}(y), \m{I}_X)$ (pull-back/restriction symbols
are mostly  omitted throughout the paper).

\subsection{Some basic  toric algebraic geometry}\la{lali}
In this section, we work over an arbitrary ground field.
Let $Z=Z(\Delta)$ be a toric variety and 
let $f: X \to Y$ be a toric map. We simply offer a list of the properties we need;
for proofs and/or references, see \ci{dMM}.

\begin{num}\la{n1}
{\rm 
{\bf  Toric affine open cover, orbit closures and partial order.}
$Z$ is covered by the open  affine
toric subvarieties $U_\zeta= \coprod_{\rho \leq \zeta} O(\rho),$ where $\zeta \in \Delta.$
We have:
$V(\zeta) := \ov{O(\zeta)} =\coprod_{\rho \geq \zeta} O(\rho),$ 
and 
 $\tau \leq \s$ iff $V(\tau) \supseteq O(\zeta).$
}
\end{num}

\begin{num}\la{n2}
{\rm 
{\bf  Toric varieties of contractible type.}
The toric variety $Z$ is said to be of contractible type if it is of the form
$(Z, z)= (U_\zeta, z_\zeta),$
where the cone $\zeta$ spans $N_{\real}.$ 
In this case, $z$ is the unique torus fixed point.
}
\end{num}

\begin{num}\la{n3}
{\rm 
{\bf  The local product structure of $Z$ along orbits.}
 Given $\zeta \in \Delta,$ any splitting   $ N \cong N_{\zeta} \oplus
N (\zeta)$  of lattices determines  a splitting   $T \cong T_\zeta \times T(\zeta)$
of tori,
and an equivariant isomorphism of toric varieties:
\beq\la{lps}
U_\zeta  \cong U_{\zeta'} \times O(\zeta),
\eeq
where $\zeta'$ is the cone $\zeta,$ viewed in $N_{\zeta} \subseteq
N_.$  One virtue of (\ref{lps}) is that $U_{\z'}$ is of contractible type and the product
assertion is useful in the context of inductive arguments; the same is true for (\ref{3d}) below.
The isomorphism (\ref{lps}) depends on the choices. The fan in $N_{\zeta'}:= N_{\zeta} \subseteq N$ given
by $\zeta'$ and its faces yields a  canonical closed embedding $(U_{\zeta'}, z_{\zeta'}) \to (U_\zeta, z_\zeta),$
 compatible with the non canonical (\ref{lps}).
}
\end{num}

\begin{num}\la{n4}
{\rm 
{\bf The local product structure of a proper toric fibration over the $U_\s$.}
Let $f$ be a  proper toric fibration. 
Let $\s \in \Delta_Y.$   
There is   a  non canonical
equivariant splitting as in (\ref{lps}),
and a  non canonical
equivariant isomorphism
of toric maps, compatible with (\ref{lps}): 
\beq\la{3d}
(f^{-1} (U_\s) \to U_\s) \cong \left( f^{-1} (U_\s') \times
O(\s) \stackrel{f_{\s'} \times {\rm Id}}\lorw  U_{\s'} \times O(\s)\right).
\eeq
The resulting natural restriction-over-$U_{\s'}$-map
$f_{\s'}$ is a toric fibration onto a base of contractible type, and  we have a natural identification
$f^{-1} (y_\s) = f_{\s'}^{-1} (y_{\s'}).$
In particular, we get a $(T_X \to T_Y(\s))$-equivariant non canonical decomposition:
\beq\la{3d1}
f^{-1} (O(\s)) \cong f^{-1} (y_\s) \times O(\s).
\eeq

}
\end{num}
\begin{num}\la{n5}
{\rm 
{\bf Canonical factorization of induced maps between orbits.}
Let  $\xi \in \Delta_X$  and consider  the natural map of tori 
$\phi: (O(\xi), x_\xi) \to (O(\ov{\xi}), y_{\ov{\xi}})$ induced by $f.$
 The image  is a  closed subtorus $i: (O'({\xi}), y_{\ov{\xi}}) \to
 (O(\ov{\xi}), y_{\ov{\xi}})),$ and  there is the following canonical
 factorization into maps of tori:
\beq\la{okoko}
\phi: O(\xi) \stackrel{a}\lorw A  \stackrel{b}\lorw B \stackrel{c}\lorw  O'({\xi}) \stackrel{i}\lorw O(\ov{\xi}),
\eeq
where: $a$ is a  toric fibration (non canonically a product projection);
$b$ is a geometric quotient map, \'etale and Galois,  by the action of
a finite Abelian subgroup of the torus   $A;$
$i$ is the natural closed embedding above; $c$ is a universal homeomorphism.
}
\end{num}

\begin{num}\la{n6}
{\rm 
{\bf The Stein factorization of a toric map.}
Let $f: X \to Y$ be a proper  toric map. 
 There is the  canonical toric  Stein factorization:
\beq\la{tsft}
f: X \stackrel{g}\to Z \stackrel{h}\to Y,
\eeq
where $g$ is a  proper toric fibration ($g_* \m{O}_X = \m{O}_Z$; surjective
with connected fibers),
and $h$ is toric finite.  The normalization of  the image $f(X)$ is 
a toric variety.
}
\end{num}

\begin{num}\la{n9}
{\rm 
{\bf Toric resolutions,  toric Chow envelopes, toric completions.}
There is a proper birational toric map $g: W \to X$ such that
$W$ is nonsingular; 
one can choose $W$ to be quasi projective, so that
$g$ is then projective.
In particular, if $f$ is proper, then there is a projective toric map $g: W \to X$
such that $h: = f\circ g$ is projective toric.

\n
There  is a toric  completion of $X,$ i.e.  an open immersion $j: X \to \ov{X}$ such that $j$ is a toric  map and $\ov{X}$ is toric complete.

\n
There are  toric completions 
$X\subseteq \ov{X}$ and $Y\subseteq \ov{Y}$ and  a proper toric map $\ov{f} : \ov{X} \to \ov{Y}$ extending $f.$  
}
\end{num}

\begin{num}\la{n7}
{\rm 
{\bf Equivariant complexes.} In what follows, we work over a field that is either finite
or algebraically closed.
Let ${\rm act}: Z \times T \to Z$ be the torus action on a toric variety. 
A complex $C$ on $Z$ is equivariant\footnote{the standard definition of equivariance
requires the usual cocycle condition; we do not need it here} if there is an isomorphism
${\rm act}^* C \cong \pi_Z^* C.$ 
We can extend this notion to the torus invariant subschemes of $Z.$
The intersection complex $\m{I}_{V(\s)}$ of the closure of an orbit is equivariant.
Given $\nu \in N,$ we have the associated co-character
$\lambda_\nu: \bb{G}_m \to T,$ and the associated notions of $\lambda_\nu$-equivariance.
If $f: X\to Y$ is a proper  toric fibration, then $Rf_*$ preserve equivariance.
}
\end{num}

\begin{num}\la{n8}
{\rm 
{\bf The retraction lemma.} This is where the notion of toric variety of contractible type
starts playing a role.
 Let $ f: X\to Y$ be a proper toric fibration onto $(Y,y)$ of contractible type
and let $C$ be an equivariant complex on $X$ (the ground field
is algebraically closed, or finite).
The natural graded  map below is an isomorphism:
\beq\la{lartzl}
 H^*(X, C) = H^*(Y, Rf_* C) \stackrel{=}\lorw (R^*f_* C)_y,
 \eeq
 where it is understood that if 
 the ground field is  finite, then we have passed
 to an algebraic closure, and we have
 an isomorphism of $G$-modules.
Special case: $C=\m{I}_X$; then,  by  coupling with  proper base change:
  \beq\la{lartzl1}
  I\!H^*(X) = (R^* f_* \,\m{I}_X)_y = H^* \left( f^{-1} (y) , \m{I}_X \right).
  \eeq
 Special case of the special case: $f= {\rm Id}_Y$;  then:
 \beq\la{spid}
 I\!H^*(Y) = \m{H}^*(\m{I}_Y)_y.
 \eeq

}
\end{num}

\subsection{The decomposition theorem for proper
toric maps over finite fields}\la{smrz}
Let $\z{f}: \z{X}\to \z{Y}$ be a proper toric map
over a finite field $\z{\bb{F}}.$   
Let $\z{f}= \z{h} \circ \z{g}:
\z{X} \to \z{Z} \to \z{Y}$ be the Stein factorization. 
For every $\zeta \in \Delta_{\z{Z}},$ define, recalling (\ref{okoko}):
\[{\rm Ev}_\zeta:= \left\{ b \in \zed \, | \, b  + \dim X - \dim V(\zeta) \; \mbox{even} \right\},
\quad \beta_\zeta:= \frac{b  + \dim X - \dim V(\zeta)}{2},
\]
\[\z{O'({\zeta})} := \z{h} ( \z{O(\zeta)}), \qquad 
\z{L}_\zeta = \z{h}_* {\oql}_{\z{O (\zeta)}} \;\;
\mbox{a lisse sheaf on}\;\;  \z{O'({\zeta})} \subseteq \z{O} (\overline{\zeta}).
\]

\begin{tm}\la{tmb} 
{\rm ({\bf DT for proper toric maps over finite fields})}

\ben
\item[\rm a)]\la{cc} 
Let $\z{f}: \z{X}\to \z{Y}$ be a proper toric map. 
There is a DT isomorphism in $D^b_m (\z{Y}, \oql):$
\beq\la{tmb4}
Rf_*  I_{\z{X}} \; \cong\;  \bigoplus_{\zeta \in \Delta_{\z{Z}}} \bigoplus_{b \in 
{\rm Ev}_\zeta} 
{I}_{\ov{\z{O'({\zeta})}}}^{s_{\zeta,b}}  (\z{L_\zeta}) (- \beta_\zeta) [-b], 
\eeq
where:    $\z{O'({\zeta})} := \z{h} ( \z{O(\zeta)});$
the sheaves $\z{L}_\zeta = \z{h}_* {\oql}_{\z{O (\zeta)}}$ 
on  $\z{O'({\zeta})}$ are lisse, semisimple,
pure of weight zero;
  the  $s_{\zeta,b} \in \zed^{\geq 0}$ are  subject to:
  
i) $s_{\zeta,b} = s_{\zeta, -b},$ for every $b\in {\rm Ev}_\zeta;$ 

ii) if $\z{f}$
is projective, then $s_{\zeta,b} \geq \sum_{l \geq 1} s_{\zeta, b+2l},$ for every $ b \geq 0$
in ${\rm Ev}_\zeta.$  

\item[\rm b)]\la{dd}
In particular: the  pure weigth zero 
$R \,\z{f_*} \,\m{I}_{\z{X}}$ is   punctually pure, even and Tate;
for every $y \in \z{Y} (\bb{F}),$ the $G$-module $(R^* f_* \m{I}_X)_y = H^*(f^{-1}(y),  \m{I}_X)$ is  pure,
even and Tate.

\item[\rm c)]\la{bb}
Let $\z{f}$ is a proper toric fibration and, for $\s \in \om_Y,$
let ${\rm Ev}_\s,$ and $\beta_\s$ as above.
There is a DT isomorphism  in $D^b_m (\z{Y}, \oql):$
\beq\la{tmb5}
Rf_*  I_{\z{X}} \; \cong\;  \bigoplus_{\s \in \om_{\z{Y}}} \bigoplus_{b \in Ev_\s} 
I_{\z{V(\s)}}^{s_{\s,b}}  (- \beta_\s) \,[-b], 
\eeq
 where   the  $s_{\s,b} \in \zed^{\geq 0}$ 
are  subject to the  conditions analogous to i) and ii)  above.
  
\een
\end{tm}

\begin{rmk}\la{supr}
{\rm 
Even though statement a) implies statements b) and c),
we prove the three assertions in the
following  order: 
 assertion b) is proved by Corollary \ref{r55}; assertion  c) is proved by Theorem \ref{dttor};
assertion   a) is proved by Theorem \ref{ptm11}, which builds on assertion c). 
}
\end{rmk}

\begin{rmk}\la{shiftz}
{\rm
We may re-write the DT isomorphisms
using the shifted $\m{I}_X$ by   setting $2k= 2\beta,$ etc.  In the case of a proper
toric fibration we get, with $s_{\s,b}$ as in (\ref{tmb4}):
\beq\la{u1}
R \, {\z{f}_*} \, \m{I}_{\z{X}} \, \cong \,  \bigoplus_{\s \in \Delta_{\z{Y}}}
\bigoplus_{k \in \zed}  
\m{I}_{\z{V(\s)}}^{ s_{\s, 2k - \dim \z{X} + \dim\z{ V(\s) }  }} (-k)
[-2k].
\eeq
While  (\ref{tmb5}) emphasizes
duality, (\ref{u1}) emphasizes  eveness.}
\end{rmk}

\begin{rmk}\la{tmc}
{\rm ({\bf DT for proper toric maps over algebraically closed fields})
Theorem \ref{tmb} over finite fields implies the analogous results
over any algebraically closed field (remove $\z{-}$ and $(-\beta)$).
In fact, assume that ${\rm char}\, K  \neq 0$  and form the tower
of field extensions $\z{\bb{F}} \subseteq \bb{F} \subseteq K$ given
by the prime subfield of $K$ and by its algebraic closure in $K.$ 
Theorem \ref{tmb}  over a finite field implies
immediately  the desired conclusions for  $K=\bb{F}.$
One then  pulls-back  further  to $K$  and concludes
when ${\rm char}\, K \neq 0.$

\n
According to  \ci{BBD}, \S6, especially \S6.1.10, Lemme 6.2.6 and 
Th\'eor\`eme 6.2.5, relating the situation over an algebraically closed field of characteristic zero to the one over an algebraic closure of a finite field,
the desired conclusion in characteristic zero
follows from the one in positive characteristic.

\n
If the ground field is $\comp$ and we use the classical topology,
then one can reach analogous conclusions by using  the theory of mixed Hodge modules.
}
\end{rmk}

\section{Decomposition theorem for proper toric maps over  finite  fields}\la{rv66}

\subsection{General  DT package over finite fields}\la{mn}
The following is surely well-known and follows easily from some of the results in
\ci{BBD}. We could not find an adequate explicit reference.

\begin{pr}\la{wdte3}
{\rm {\bf (DT and RHL over $\z{\bb{F}}$)}} Let $\z{f}: \z{X} \to \z{Y}$ be a proper map of  separated $\z{\bb{F}}$-schemes of finite type, let $\z{P}$ be a pure perverse sheaf
of weight $w$ on $\z{X}$. Then:  the   direct image complex $R\, \z{f}_* \z{P}$ is pure
of weight $w$
and splits non canonically into the direct sum $\bigoplus_b {^p\!H}^b(R \z{f}_* \z{P}) [-b]$
of its shifted perverse direct image  complexes; the perverse sheaves ${^p\!H}^b(R
\, \z{f}_* \z{P})$ 
are pure of weight $w+b$ and admit the  canonical decomposition by supports as the direct sum $\bigoplus_{\z{\frak Y}} 
IC_{\z{\frak Y}} (\z{L}_{b,\z{\frak Y}})$  of finitely many intersection complexes
of  $\z{\bb{F}}$-integral subvarieties  $\z{\frak Y}\subseteq \z{Y}$, with coefficients  pure lisse
sheaves $\z{L}_{b,\z{\frak Y}}$ of weight $w+b- d_{\z{\frak Y}}$ on suitable Zariski dense open
subsetes $\z{\frak Y}^o \subseteq \z{\frak Y}$. 

\n
Assume, in addition, that $\z{f}$ is projective and let
$\z{\eta}$  be the first Chern class of an $\z{f}$-ample line bundle on $\z{X}$.
 Then the relative hard Lefschetz theorem (RHL) holds: for every $i \geq 0,$
  the cup product map
 $ \z{\eta}^i: {^p\!H}^{-i} (R\,\z{f}_* \z{P}) \lorw
 {^p\!H}^{i} (R\,\z{f}_* \z{P})(i) $ on the perverse cohomology sheaves of the direct image  is an isomorphism.
 \end{pr}
{\em Proof.} By Deligne's fundamental result \ci{BBD}, 5.1.14, the direct image $R\, \z{f}_* \z{P}$ is pure  of weight $w.$ The perverse sheaves ${^p\!H}^b(R \z{f}_* \z{P})$ are pure of weight $w+b$ by \ci{BBD}, Thm. 5.4.1.
They split as indicated by virtue of \ci{BBD}, Cor. 5.3.11, coupled with a straightforward
Noetherian induction.

\n
In the projective case, the RHL is  \ci{BBD}, Theorem 5.4.10, and the splitting
into the direct sum of shifted perverse cohomology sheaves is a formal consequence
of RHL via the Deligne-Lefschetz criterion.

\n
The DT for a proper map can be derived formally as follows.

\n
By using  \ci{BBD}, Cor. 5.3.11, we first observe that we may assume 
that $\z{P}$ is  the intermediate extension back to $\z{X}$
of its own restriction to any Zariski-dense open subvariety.

\n
We choose a 
Chow envelope $\xymatrix{\z{h}: \z{W} \ar[r]^{\z{g}} & \z{X} \ar[r]^{\z{f}} & \z{Y}}$
of the map $\z{f}$, so that $\z{g}$ and $\z{h}$ are projective and there is
a Zariski dense open subvariety $\z{U} \subseteq \z{X}$ over which 
$\z{g}$ is an isomorphism. 

\n
Let $\z{P'}$ be the corresponding intermediate extension
to $\z{W}$. Since $\z{g}$ is projective, we can apply the Deligne-Lefschetz 
splitting and deduce
that $\z{P}$ is a direct summand of $R\,\z{g_*} \z{P'}$.

\n
The conclusion follows 
 by applying what we have already proved for the projective
$\z{g}$ and to $\z{h}$, and then by using
the relation $R\,\z{h}_* = R\,\z{f}_* \circ R\,\z{g}_*$ to compare terms.
\blacksquare

\begin{rmk}\la{dtfp}{\rm
The lisse sheaves $\z{L}$  on $\z{\frak{Y}^o}$ become semisimple,
after pull-back to
each integral component of $\frak{Y}^o.$ As we shall see, in the case of
proper toric maps with $\z{P}:= I_{\z{X}},$ the subvarieties $\z{\frak Y},$
being torus orbits, 
are geometrically integral, and, in addition, the coefficients $\z{L}$ turn out to be are already semisimple.
See Theorems \ref{dttor}, \ref{ptm11}. 
The indecomposable pure perverse sheaves on  a $\z{Y}$ are described in
\ci{BBD}, Prop. 5.3.9]. Due to my ignorance, I  do not know if the pure lisse $\z{L}$
may admit indecomposable direct summands that present  Jordan block-type
factors of type $E_n$ with $n\geq 2:$ that would be the only obstacle to 
the semisimplicity  of the pure lisse coefficients $\z{L}$  on $\z{Y}.$
}
\end{rmk}

\subsection{First  toric consequences of Proposition \ref{wdte3}}\la{test1}

\begin{lm}\la{warg}
Let $\z{Z} \stackrel{\z{u}}\to \z{U} \stackrel{\z{j}}\to \z{V}$
be maps of schemes (separated and of finite type)  over a finite field and let $\z{C}$ be a pure complex
of weight $w$ on
$\z{V}.$  Assume that:
$\z{Z}$ is complete, $\z{j}$ is an open immersion, $u^*: H^*(U,C) \to H^*(Z,C)$
is an isomorphism, $\z{V}$ is smooth.
If  the -automatically pure-- submodule of lowest weight $w$ of $H^*(V,C)$ is  even and Tate, then so is $H^*(U,C).$
\end{lm}
{\em Proof.} This is standard; we freely use  basic weight theory (\ci{BBD}, 5.1.14). It is enough to show that $j^*$ is surjective.
We have that:
$H^*(U,C) \cong H^*(Z,C)$ has weights $\leq w,$ because $\z{Z}$ is complete,
so that the direct image coincides with the extraordinary direct image,
under which the property of having weights $\leq w$ is stable.
$H^*(U,C)$ has weights $\geq w$ because the property of having weights $\geq w$ is stable under direct image. By combining the two weight inequalities above,
we see that  $H^*(U,C)$ is pure of weight
$w.$  Let $\z{i}$ the closed embedding complementary to $\z{j}.$
In view of the fact that the property of having weight $\geq w$ is stable under extraordinary inverse image,
the long exact sequence of relative cohomology $H^*(V,C)
\to H^*(U,C) \to H^{*+1}(V, i_!i^!C)$ shows
that $j^*$ is surjective.
\blacksquare

\begin{lm}\la{t5}
Let $\z{X}$ be a toric variety over a finite field. The
intersection complex $\m{I}_{\z{X}}$ is pure of weight zero and $\om$-constructible.
\end{lm}
{\em Proof.} Both are well-known. We offer a proof of purity for toric
varieties based on the DT. Of course, purity of the intersection complex
for any Variety over a finite field
is a result of O. Gabber.
We also prove $\om$-constructibility as we need some of the details of the proof
in the proof of Lemma \ref{vbg}.

\n
Let $\z{g} : \z{W}\to \z{X}$ be a proper toric resolution of the singularities of
$\z{X}.$ Proposition \ref{wdte3} implies that $\m{I}_{\z{X}}, $ being  a direct 
summand of the pure weight zero $R \z{g}_* {\oql}_{\z{W}},$ is pure of
weight zero.

\n
The proof of $\om$-constructibility is by induction on $n:= \dim \z{X}.$
If $n=0,$ then we are done. Assume the desired conclusion holds
for every toric variety of dimension at most $n-1.$

\n
Since the open  sets of the form $\z{U_\s}$ are union of orbits
and cover $\z{X}$ as $\s$ ranges in $\Delta_X,$ we may assume that  
 $\z{X}= \z{U_\s}.$ In view of the local product structure $\z{U_\s} \cong \z{U_{\s'}}
 \times \z{O(\s)},$ we may also assume that $(\z{X}, \z{x})$ is of contractible type.
 
 \n
 Since $\z{x}$ is now an orbit, we may replace $(\z{X}, \z{x})$
 with  $\z{X} \setminus \z{x}.$ Now,
 $\z{X}$ is  covered by the affine open sets of the form $\z{U_{\tau}} \cong
 \z{U_{\tau'}} \times \z{O(\tau)} ,$ with $\tau <\s,$ so that
 $\dim \z{U_{\tau'}} \leq n-1$, and we are done.
 \blacksquare
 
 \begin{lm}\la{vbg}
 Let $\z{f}: \z{X} \to \z{Y}$ be a proper toric fibration over a finite
 field.
 The direct image complex $R\, \z{f_*} \,\m {I}_{\z{X}}$ is $\om$-constructible.
 \end{lm}
 {\em  Proof.}  As in the proof of Lemma \ref{t5}, 
we may  assume that  $\z{Y}= \z{U_\s}.$ 

\n
According to (\ref{3d}), the intersection complex
of $\z{f}^{-1} (\z{U_\s})$ is a pull-back from the factor
$\z{f}^{-1} ( \z{U_{\s'}} ).$ It follows that we may assume that 
$(\z{Y}, \z{y})$ is of contractible type.

\n
We conclude the proof by arguing by induction on the dimension
of the base of contractible type as 
 in the proof of Lemma \ref{t5}.
  \blacksquare

\begin{lm}\la{r44}
Let $\z{f}: \z{X} \to \z{Y}$ be a proper toric fibration onto a base of contractible type
$(\z{Y}, \z{y}),$ all over a finite field. Then we have natural 
isomorphisms of graded $G$-modules:
\beq\la{r440}
I\!H^*(X) = (R^* f_* \m{I}_X)_y= H^* (f^{-1} (y), \m{I}_X).
\eeq
In particular, we have the natural isomorphism of graded $G$-modules:
\beq\la{r441}
I\!H^*(Y) =  \m{H}^*(\m{I}_Y)_y.
\eeq
The graded $G$-modules above are pure, even and Tate.
\end{lm}
{\em Proof.} The first two assertions follow from the retraction 
lemma  (\ref{lartzl1}) and (\ref{spid}) and proper base change.

\n
We turn to the proof that  the $G$-modules (\ref{r440}) are pure even and Tate. 

\n
By taking a projective toric resolution of $\z{X}$ as in the proof of Lemma
\ref{t5}, we may assume that $\z{X}$ is smooth and quasi projective.
In particular, we are now dealing with cohomology, instead of intersection cohomology. 

\n
Choose a toric open embedding $j: \z{X} \to \ov{\z{X}}$ with $\ov{\z{X}}$
toric smooth and projective. 

\n
It is well known  that the graded $G$-module $H^*(\ov{\z{X}})$
is pure of weight zero, even and Tate: it is generated by algebraic cycle classes
\ci{Ful93}. 

\n
We conclude by applying Lemma \ref{warg} to $\z{f}^{-1} (\z{y}) 
\to \z{X} \to \ov{\z{X}}.$
\blacksquare

\begin{rmk}\la{fz}
{\rm ({\bf Description of $H^* (f^{-1} (y), \m{I}_X)$, $f$ proper toric)}
Let $f: X\to Y$ be a proper toric map over an algebraically closed  field and let
$y \in Y$ be a closed point.  Let $f=h \circ g: X \to Z \to Y$ be the toric Stein factorization.  The fiber $f^{-1} (y)$ is a disjoint union of the finitely many fibers
 $g^{-1}(z_l),$ with $\{z_l \} = h^{-1} (y).$ 
 Clearly, $H^*(f^{-1}(y), \m{I}_X)=
 \oplus_l H^* (g^{-1} (z_l), \m{I}_X).$ 
Fix $z_l$ and take    the $U(\zeta)$ with $z \in O(\zeta).$
 By using the local product structure  (\ref{3d}) over $U_{\zeta},$ we may assume that
 $z_l= z_\zeta$ and that  $(Z,z_\zeta)$
 is of contractible type. By combining the proofs of lemmata \ref{r44} and \ref{t5},
 we obtain the following description:
 
\n
{\em $H^* (f^{-1}(y), \m{I}_X)$ is a  finite direct sum
of subquotients   of the graded vector spaces 
  $H^*(\ov{W_l}, \oql)$ given by the cohomology of  nonsingular projective
  toric varieties $\ov{W_l}$.
  }
 \n
 The description
remains valid in the context of $G$-modules if the ground field is the algebraic closure
of a finite field.  Similarly,  over $\comp,$  in the context of the classical Euclidean topology
with  the rational mixed Hodge structures of the theory of mixed Hodge modules.
}
\end{rmk}

\begin{cor}\la{r55}
Let $\z{f} : \z{X} \to \z{Y}$ be a proper  toric map  over a finite field. Then
the pure weight zero  complexes $\m{I}_{\z{X}}$ and $R \, \z{f_*} \,
\m{I}_{\z{X}}$  are  punctually pure,
even and Tate.
\end{cor}
{\em Proof.}  The second statement implies the first one by taking
$\z{f} = {\rm Id}_{\z{X}}.$

\n
Since the desired conclusions are stable under direct images via finite toric maps,
in view of the toric Stein factorization of the map $\z{f},$ we may assume
that $\z{f}$ is a proper toric fibration.

\n
Since the desired conclusions are statements about the stalks
and the direct image complex is $\om$-constructible, it is enough
to verify that for every distinguished point $\z{y_\s}$ of any orbit 
$\z{O(\s)} \subseteq \z{Y},$ the graded $G$-module
$(R^* f _* \m{I}_X)_{y_\s}$ is pure, even and Tate. This follows
immediately from Lemma \ref{r44}.
\blacksquare

\begin{rmk}{\rm The punctual purity etc. of the intersection complex of a  toric variety
is proved in \ci{DL}. The statement for the direct image seems new.}
\end{rmk}

\subsection{Proof of Theorem \ref{tmb}.c on proper  toric fibrations}\la{pfdttr}
The following theorem proves Theorem \ref{tmb}.c.

\begin{tm}\la{dttor} {\rm ({\bf Proper toric fibrations: DT  semisimplicity and RHL over $\z{\bb{F}}$)} }

\n
Let $\z{f}: \z{X}\to \z{Y}$ be a proper  toric fibration.
There are    isomorphisms in $D^b_m (\z{Y}, \oql):$
\beq\la{fgoo}
R\, \z{f}_* I_{\z{X}} \cong \bigoplus_ {\s \in \Delta_{\z{Y}} } 
\bigoplus_{b \in \zed} \, \z{J}_{\s,b} \,[-b],
\eeq
\beq\la{333gbh}
\z{J}_{\s,b} \cong 
\left\{
\begin{array}{ll}
0&\mbox{$b  + \dim X - \dim V(\s)$ odd,}\\
\z{I}_{\s}^{s_{\s,b}} 
\left(- \frac{b  + \dim X - \dim V(\s)}{2} \right) &\mbox{$b  + \dim X - \dim V(\s)$ even.}
\end{array}
\right.
\eeq

\n
The  integers $s_{\s,b}$ are subject to the following relations:
\ben
\item[i)]
 $s_{\s,b} = s_{\s,-b}$ (Poincar\'e-Verdier duality); if
 $b+\dim{X} - \dim V(\s) $ is odd,  then $s_{\s,b} =0.$
 
 \item[ii)]
 If $\z{f}$ is projective, then
  for every $b\geq 0$, $s_{\s,b} \geq \sum_{l\geq 1} s_{\s, b+2l}$ (RHL).
   \een
 \end{tm}

\n
{\em Proof.}
Proposition \ref{wdte3}, coupled with the $\om$-constructibility
of the direct image complex  in Lemma \ref{t5}, implies the existence of an isomorphism
(\ref{fgoo}), where each  $\z{J_{b,\s}}$   is of the form
$\z{I_{V(\s)}}(\z{L_{b,\s}})$, with $\z{L_{b,\s}}$ lisse,  pure of weight $b+ \dim X  
-\dim V (\s)$
on $\z{O(\s)}$. 

\n
By $\om$-constructibility, for each $j,\s$,  the mixed
sheaf $\z{R^j_\s}$ (see the end of \S\ref{notz} on  notation)  is lisse   on the orbit $\z{O(\s)}.$

\n
The local product structure  of $\z{f}$  (\ref{3d1}) implies that 
$\z{R^j_\s}$ is given by the  continuous representation:
\beq\la{rt51}
\rho^j_\s: \pi_1 (\z{O(\s)}, y_\s) \lorw  \pi_1 (\z{y_\s}, y_\s) \lorw {\rm GL}_{\rat_\ell}
\left( R^j_{y_\s} \right),
\eeq
where the first homomorphism stems from the natural constant  $\z{\bb{F}}$-map 
$\z{O(\s)} \to  \z{y_\s},$ 
and the second one from the $G$-module structure on $R^j_{y_\s}:$ this means that
the lisse sheaf  $\z{R^j_\s}$ is the pull-back of a lisse sheaf on $\z{y_\s},$ namely,
the $G$-module $R^j_{y_\s}.$

\n
By combining with Corollary \ref{r55}, 
we see that
$\z{R^j_\s}$ is zero for $j$ odd and that, for $j$ even, we have
$\z{R^j_\s} \cong {\oql^r}_{\z{O(\s)}} (-j/2),$ for some $r\geq 0.$

\n
The already-proved existence of an isomorphism (\ref{fgoo})
as above,  implies that $\z{L_{b,\s}}$ is a direct summand
of $\z{R^j_\s}$ for $j = b+  \dim X  - \dim V(\s )$, so that we have
$\z{L_{b,\s}} \cong {\oql^s}_{\z{O_\s}} (-j/2),$  for some integer $s,$ equal to zero
if $j$ is odd. This proves (\ref{333gbh}),  as well as  the second part of assertion i).

\n
Finally, given (\ref{fgoo}) and (\ref{333gbh}), assertion ii) and the first part of i)  are immediate consequences of  the relative hard Lefschetz theorem,
and of  Verdier duality,  respectively.
 \blacksquare

\subsection{Proof of Theorem \ref{tmb}.a on proper toric maps}\la{ptm}
Let  $\z{f}: {\z{X}} \to {\z{Y}}$ be a proper  toric map  over a  finite field $\z{\bb{F}}$ and let
$\z{f}= \z{h} \circ \z{g}$ its  Stein factorization.
Since $\z{g}$ is a proper toric fibration, i.e. the subject of Theorem \ref{dttor},
 we turn our attention to  the finite  map ${\z{h}}: {\z{Z}} \to {\z{Y}},$ with
 image  $\z{f}(\z{X}),$  a closed subvariety of $\z{Y}.$ 

 Given $\zeta \in \Delta_{\z{Z}},$ we have the map of tori ${\z{h(\zeta)}}:
{\z{O(\zeta) }} \to {O({\ov{\zeta})}},$
with image $\z{O' ({{\zeta}}) }\subseteq  O({\ov{\zeta}}).$  We denote the
evident resulting finite map on the closures by 
$\ov{\z{h(\zeta)}}: {\ov{\z{O(\zeta)}}}
\to {\ov{\z{  O'({{\zeta}})  }}}.$

\begin{lm}\la{rt5m}
The ${\oql}$-adic sheaf $\z{L_{\zeta}}:= {\z{h(\zeta)}}_* {\oql}_{\z{O(\zeta)}}$ on
$\z{O' ({{\zeta}}) }$
 is lisse, semisimple,  pure, punctually pure of weight zero;  the geometric monodromy has eigenvalues roots of unity.
 
 \n
 The direct image 
 $R \, {\z{h_*}} \,\z{\m{I}_{V(\zeta)}}=  \m{I}_{{\z{O' ({{\zeta}})  }}} 
 ( \z{L_{\zeta}} )$ 
 is pure of weight zero, punctually pure, even and Tate. 
\end{lm}

 \n
{\em Proof.}  The map of tori $\z{h(\zeta)}$ admits 
the canonical factorization (\ref{okoko}).
Recalling that the map being factored is finite: the map $\z{a}$ must be the identity.
 The map $\z{b}$ is  a quotient
by a finite abelian group, $\Gamma:= \ke \, \z{b},$ whose order is not divisible
by ${\rm char}\, \z{\bb{F}};$ in particular it is \'etale,
and Galois.  The map $\z{c}$ is a universal homeomorphism.
 The map $\z{i}$ is the evident  closed
embedding.

\n
Since $\z{a}, \z{c}$ and $\z{i}$ are universal homeomorphisms onto their image,
they do not effect the direct image calculations and can be ignored.

\n
It follows that $\z{L_\zeta}$ is naturally identified the ordinary direct
image sheaf $\z{b}_* \oql$, with $\z{b}$ the \`etale quotient by the  action
of the finite abelian group $\Gamma.$
All the listed properties of $\z{L_\zeta}$  follows easily from this description.

\n
The equality  statement about the direct image complex  follows
by observing that we have the following natural identifications:
\beq\la{dkfj}
R \, {\z{h_*}} \,   \m{I}_{\z{V(\zeta)}} = \z{h_*} \,  \m{I}_{\z{V(\zeta)}}  = 
 \ov{\z{h (\zeta)}}_* \m{I}_{ \z{V(\zeta)}  }
 =
 \m{I}_{\ov{\z{O' ({{\zeta}})  }}} 
 ( \z{L_{\zeta}} ):
 \eeq
the first one is because $h$ is finite. The  second one is
because the third term is merely
a re-writing of the second; the third identification 
follows from the  fact that the direct image under
the finite map $\ov{ \z{h(\zeta)}}$  is $t$-exact for the middle perversity $t$-structure and hence
preserves intersection complexes with twisted coefficients; clearly, the coefficients
of the direct image can be read on a Zariski-dense open subset, so that they
are given by $\z{L_\zeta}.$ 

\n
Finally, the last statement follows by the just-established equality
and from Corollary \ref{r55}, applied to $\z{X}:= \z{V(\zeta)}: $ 
for the  properties in question are stable under finite direct image.
\blacksquare

\begin{rmk}\la{rmkchar}
{\rm 
The rank of the local system $\z{L_\zeta}$ is the cardinality of
the abelian group $\Gamma,$ which depends on the characteristic of the finite field
$\z{\bb{F}}.$
}
\end{rmk}

\smallskip
We can now show that 
Theorem \ref{dttor} for proper toric fibrations
over $\z{\bb{F}}$  has  the following natural  counterpart  for proper toric maps over $\z{\bb{F}},$ which, in turn, establishes Theorem \ref{tmb}.a.

\begin{tm}\la{ptm11}
{\rm ({\bf Proper toric maps: DT semisimplicity and RHL over $\z{\bb{F}}$)} }
Let ${\z{f}}: {\z{X}} \to {\z{Y}}$ be a proper toric map over the finite field $\z{\bb{F}}.$
There are  isomorphisms in $D^b_m (\z{Y}, \oql):$
\beq\la{po}
R \, \z{f_*} I_{\z{X}} \cong  \bigoplus_{\zeta \in \Delta_{\z{Z}}} \bigoplus_{b \in \zed}
\ov{\z{J_{\zeta,b}}} [-b], 
\eeq
\beq\la{3gh}
\ov{\z{J}_{\zeta,b}} \cong 
\left\{
\begin{array}{ll}
0,&\mbox{$b+ \dim X - \dim V(\zeta)$ odd,}\\
\z{IC_{\ov{\z{O'({{\zeta} )}}}}^{s_{\zeta,b}}} (\z{L_{\zeta}})
\left(- \frac{b+ \dim X  - \dim  V(\zeta)}{2} \right), &\mbox{$b+
\dim X - \dim V(\zeta)$ even.}
\end{array}
\right.
\eeq
The complex $R\, \z{f}_* I_{\z{X}}$ is pure of weight zero,  punctually pure,
even and Tate.

\n
The  integers $s_{\zeta,b}$ are subject to the following relations:
\ben
\item[i)]
 $s_{\zeta,b} = s_{\zeta,-b}$ (Poincar\'e-Verdier duality);
 if $b+ \dim{Z} - \dim  V(\zeta)$ is odd, then $s_{\zeta,b} =0.$
 
 \item[ii)]
 If $\z{f}$ is projective, then
  for every $b\geq 0$, $s_{\zeta,b} \geq \sum_{l\geq 1} s_{\zeta, b+2l, \zeta}$ (RHL).
   \een
\end{tm}
{\em Proof.} 
We have the Stein factorization $\z{f}= \z{h} \circ \z{g}.$
We first apply   Theorem
\ref{dttor}.(\ref{333gbh}) to  the proper fibration $\z{g}.$ We
form  the $R \, \z{g}_*$ of each resulting direct summand. We apply  Lemma 
\ref{rt5m}, which remains   valid also after  arbitrary Tate twists, to each resulting term. \blacksquare

Authors' addresses:

\smallskip
\n
Mark Andrea A. de Cataldo,
Department of Mathematics,
Stony Brook University,
Stony Brook,  NY 11794, USA. \quad 
e-mail: {\em mark.decataldo@stonybrook.edu}

\end{document}